\documentclass[11pt]{article}
 \usepackage{amssymb,amsmath}
 \usepackage{amsthm}
 \newtheorem{theorem}{Theorem}[section]
 \newtheorem{theorem*}{Theorem A\!\!}
 \newtheorem{proposition}{Proposition}[section]
 \newtheorem{corollary}{Corollary}[section]
 \newtheorem{proposition*}{Proposition A\!\!}
 \newtheorem{corollary*}{Corollary A\!\!}
 \newtheorem{lemma}{Lemma}[section]

\begin{document}

 
 \title{Analysis on flag manifolds and Sobolev inequalities}

 \author{Bent \O rsted}
\date{}
 \maketitle
 \begin{abstract}
\centerline {\it To Joseph A. Wolf, with admiration} 
\vspace{2 mm} \ \\
Analysis on flag manifolds $G/P$ has connections to both representation
 theory and geometry; in this paper we show how one may derive some new
 Sobolev inequalities on spheres by combining rearrangement inequalities 
 with analysis of principal series representations of rank-one
 semisimple Lie groups. In particular the Sobolev inequalities
obtained involve hypoelliptic differential operators as 
opposed to elliptic ones in the usual case. One may hope that
 these ideas might in some
form be extended to other parabolic geometries as 
well.\footnote{ 
\newline Key
words: Sobolev inequalities, parabolic geometry, principal series representations }
\end{abstract}
 \section* {Introduction}

J. A. Wolf has worked in and has made lasting contributions to large areas of
mathematics, including Riemannian geometry, complex geometry,
representations of Lie groups, infinite-dimensional
Lie groups, and the role of flag manifolds from
many points of view. Of particular importance is his study of boundary
components of Riemannian Hermitian symmetric spaces where the role
of parabolic subgroups $P$ in semisimple Lie groups $G$ is elucidated.
At the same time he has treated many aspects of induced representations
${\rm Ind}_P^G(W)$ and the relation between the geometry of
the flag manifold 
$S = G/P$ and the analysis of representations in vector bundles over $S$.

In this paper we shall consider such parabolically induced representations 
$ \pi_{\lambda} = {\rm Ind}_P^G(L_{\lambda})$ where $\lambda$ is a parameter
for a line bundle over $S$. The aim is to understand the relation between the
detailed structure of $\pi_{\lambda}$, in particular the restriction to
a maximal compact subgroup $K \subset G$, the eigenvalues
of some standard intertwining operators, and certain Sobolev inequalities
in the space of sections of the corresponding line bundle involving natural
differential operators. These operators will reflect the geometry of $S$, in
particular understood as a {\it parabolic geometry}, such as for example
conformal geometry or the usual CR-geometry; in addition, there will be a quaternionic
analogue of the usual CR-geometry as well as an octonionic analogue.
In particular the operators that arise are hypoelliptic differential
operators as opposed to elliptic ones in the usual case. 

Some of our results will be known to experts, e.g. in parabolic geometry
with regards to covariant differential operators, 
and in representation theory with connections
to the structure of principal
and complementary series representations,
but probably their combination is new,
in particular the ensuing Sobolev inequalities. Prominent examples
of the differential operators in question will be the Yamabe operator
appearing in conformal differential geometry, and also the
CR-Yamabe operator from classical CR-geometry \cite{jl2}.  

The main result is Theorem 2.1, which gives a bound on the entropy of
a function on a sphere, viewed as a flag manifold for a rank-one simple
Lie group $G$, in terms of the smoothness of the function; the point is here
that the smoothness is only measured in certain directions in each tangent space,
corresponding to a natural distribution. This distribution is the structure that
is directly related to the structure of $G$, and it provides the relevant
parabolic geometry of the sphere in question.
    
Also, we have included in
the final section a new proof of the logarithmic Sobolev inequality
by L. Gross \cite{gr1} for the Gauss measure; this proof has the advantage of
potentially
extending to a similar inequality on the Heisenberg group, following
as a corollary to our main Theorem 2.1.

{\bf Dedication.} It is a pleasure to let this paper be part of a tribute to
Joseph A. Wolf for his mathematical work and continued energy in
revealing new insights, for his contributions as teacher and colleague  
to differential and complex geometry, Lie groups,
representations, and knowledge in general. 
                 
 \section{Geometry of the rank-one principal series}

Let $G$ be a noncompact connected semisimple Lie group with finite center;
later we shall assume that $G$ is of split rank-one. The Lie algebra of $G$ is
denoted by $\mathfrak g$.
$K = G^{\theta} \subset G$ is a maximal compact subgroup corresponding to the
Cartan involution $\theta$, and we use the same letter for the differential
giving rise to the decomposition
$$\mathfrak g = \mathfrak k \oplus \mathfrak s$$
into $\pm 1$ eigenspaces respectively. We fix a maximal abelian subspace
$\mathfrak a \subset \mathfrak s$ and for $\alpha \in \mathfrak a^*$
let $$ \mathfrak g_{\alpha} = \{X \in \mathfrak g \,|\,(\forall H \in \mathfrak a) [H,X] = \alpha(H)X \},$$
so we have the set of roots $\Delta = \{\alpha \in \mathfrak a^* \setminus \{0\}\,|\, \mathfrak g_{\alpha}
\neq \{0\} \}.$ We choose a positive system $\Delta^+ \subset \Delta$.
As usual \cite{wal} \cite{w} we have the spaces $\mathfrak m \oplus \mathfrak a =
\mathfrak g_0$ and $\mathfrak n = \bigoplus_{\alpha \in \Delta^+} \mathfrak g_{\alpha}$
as well as the corresponding analytic subgroups
$A = {\rm exp}\,\mathfrak a$, $N = {\rm exp}\,\mathfrak n$, and the
minimal parabolic subgroup $P = MAN$, where $M = Z_K(\mathfrak a)$,
the centralizer of $\mathfrak a$ in $K$,  
has Lie algebra $\mathfrak m$.

We shall be interested in representations induced from characters of $P$,
({\it scalar principal series representations} of $G$)
namely with $2 \rho = \sum_{\alpha \in \Delta^+} m_{\alpha} \alpha$,
$m_{\alpha} = {\rm dim}\,\mathfrak g_{\alpha}$, and $\lambda \in \mathfrak a_{\mathbb C}^*$
we consider $V_{\lambda} = {\rm Ind}_P^G(\chi_{\lambda})$, the representation
space of sections of the line bundle over $S$ corresponding to the
character
$$\chi_{\lambda}(man) = a^{\rho + \lambda}.$$ 
The action of $g \in G$ is by left translation and denoted by $\pi_{\lambda}(g)$
and we sometimes also use the name $\pi_{\lambda}$ for the representation.  
We identify $S = G/P = K/M$ and realize our induced representation
in $L^2(S)$ (normalized $K$ - invariant measure) as
$$\pi_{\lambda}(g)f(\xi) = a(g^{-1}\xi)^{-\lambda - \rho} f(g^{-1} \cdot \xi),$$
where $\xi = k M \in S$, $g\cdot \xi$ denotes the $G$-action on $S$, and
$a(g \xi)$ denotes the \\$A$ - component in the $KMAN$ decomposition
of $gk$. For $\lambda \in i\mathfrak a^*\, ,\pi_{\lambda}$ is unitary;
and the smooth vectors are just the smooth functions on $S$.          

Our aim is to combine some of the results in \cite{boo} and \cite{jw} with estimates
by E. Lieb \cite{l} and W. Beckner \cite{be1}, \cite{be2} in order to
obtain some new Sobolev type inequalities; they rely on some
classial rearrangement inequalities (see the Appendix), and may be thought of as
new instances of logarithmic Sobolev inequalities as found and
studied in particular by L. Gross \cite{gr1}, \cite{gr2}.
While we expect many of the results to hold in greater generality, 
we shall from now on consider the rank-one case, i.e. assume $A$
is one-dimensional, so that in particular the parabolic subgroup $P$
is also a maximal parabolic subgroup.

This means that we are dealing with (up to coverings) four
cases:
\begin{itemize}
\item {\it the real case} $G = SO_o(1,n+1),$
\item{\it the complex case} $G = SU(1,n+1),$
\item{\it the quaternionic case} $G = Sp(1,n+1),$
\item{\it the octonionic case} $G = F_4,$
\end{itemize} 
where the last case is the real form with $K = {\rm Spin(9)}$ (this $G$ could
be thought of as an analogue of octonionic $3 \times 3$ matrices
preserving a form of signature $(1,2)$). In these four cases the flag
manifold $S$ is a sphere of dimension $n, 2n+1, 4n+3, 8+7 = 15$
respectively.
Let ${\cal H}_k$ be the space of spherical harmonics of degree $k$,
i.e. homogeneous harmonic polynomials of degree $k$, restricted to $S$.
They form a representation of $K$, irreducible in the real case.                    
In the other cases we shall explicitly decompose ${\cal H}_k$ into
irreducible representations of $K$ in order to control the constants in
our Sobolev estimates. In all cases the representation $\pi_{\lambda}$
restricted to $K$ is the same as $L^2(S)$ and hence may
be identified with the direct sum of all the ${\cal H}_k,\, k \geq 0$.

Just like we can find the spectrum in $L^2(S)$ of $\Delta$, i.e., the usual Laplace--Beltrami
operator on $S$, so can we find the spectrum of the standard 
Knapp--Stein intertwining operators - see below for more on
intertwining operators; this is where we use \cite{boo} and \cite{jw},        
which we now recall.

The spherical principal series representations $\pi_{\lambda}$ depend on a
single parameter $\lambda$, which is in natural duality with $\pi_{-\lambda}$
via the invariant pairing given by integration over $K$:
$$<f,f^*> = \int_K f(k)f^*(k) dk = \int_S f(\xi) f^*(\xi) d\xi$$
where $f(k) = f(kM) = f(\xi)$ is a section of $\pi_{\lambda}$, resp.
$f^*$ a section of $\pi_{-\lambda}$. 

A central object in the representation theory of semisimple Lie groups
is that of an {\it intertwining operator}, meaning a $G$-morphism between
two modules (or a morphism for the action of the Lie algebra). There are
several standard constructions of such operators, and their analysis is
the key to many results about the structure of modules. 

We shall be
interested in intertwining operators both of integral operator type
and differential operator type; the latter occur typically as residues
of meromorphic families of intertwining operators of integral operator
type. For our purposes the relevant intertwining operators are
$$I_{\lambda} : V_{\lambda} \rightarrow V_{-\lambda},$$
where
$$ I_{\lambda} \pi_{\lambda}(g) = \pi_{-\lambda}(g) I_{\lambda}$$
for all $g \in G$ (or the analogue for the infinitesimal action of
the Lie algebra).       
Note that in this case the invariant pairing above gives rise to
an invariant Hermitian form on $V_{\lambda}$, namely
$$(f,f) = <f, I_{\lambda}\overline{f}>$$  
for $f \in V_{\lambda}$.
   
We choose an element $H_0 \in \mathfrak a$
with $\alpha(H_0) = 1$, where $\Delta^+ = \{\alpha\}$ (real case), or 
$\Delta^+ = \{\alpha, 2\alpha\}$ (remaining cases) where we have
$m_{2\alpha} = 1, \, 3, \, 7,$ resp. in the three cases complex,
quaternionic, and octonionic. 
Hence we may
identify the parameter $\lambda + \rho$ with a (in general) complex number $\nu$;
this is done by setting $(\lambda + \rho)(H_0) = \nu$.
We call the corresponding representation space $Y_{\nu}$, which may be
identified with $L^2(S)$ as a representation of $K$.

The relevant geometry now is the CR-structure on the sphere $S$.
Let us first in some detail recall the usual CR-structure on $S^{2n+1}$
and the relation to the Heisenberg group $H^{2n+1}$ \cite{jl}; this is the complex
case in our list above. 

We parametrize the Heisenberg group as $\mathbb C^n \times \mathbb R$ with  
the group product
$$(z,t)(z',t') = (z+z',t+t' + 2 {\rm Im}\,z \cdot \overline{z}'),$$
where $z \cdot \overline{z}' = \sum_{j=1}^nz_j\overline{z'_j}$.
We parametrize the Lie algebra in the same way and consider the
horizontal subspace $H_0$ defined by $t = 0$. By left translation $L_g, \, g \in H^{2n+1}$, 
the distribution $H_g = dL_g(H_0)$ on $H^{2n+1}$ is defined corresponding to the 
CR-structure, and the corresponding CR-Laplacian is
$$\Delta_b = \Delta_x + \Delta_y + 4(y\cdot  \nabla_x - x \cdot \nabla_y)\frac{\partial}{\partial t}
+ 4 z\cdot\overline{z} \frac{\partial^2}{\partial t^2}$$       
in terms of the usual Laplacian and gradient in the variables $x,y \in \mathbb R^n, \, z = x + iy$.
This operator is hypoelliptic, and it 
corresponds to the general definition in terms of a contact form; here
$\theta = dt + \sum_{j=1}^n(iz_jd\overline{z}_j - i\overline{z} dz_j)$, and
its Levi form $L_{\theta}(Z,W) = -id\theta(Z, \overline{W})$. Then
$$\Delta_b = d_b^* d_b, \, d_b = \pi \circ d: C^{\infty}(M) \rightarrow H^* = \theta^{\perp}
\subset T^*(M)$$
in general on a CR-manifold; here $M = H^{2n+1}$. Here $d$ is the
usual exterior differentiation, and $\pi$ is the dual to 
the injection $H \subset TM$. The Levi form defines the inner
product on the distribution, which is what is needed to form the adjoint of
the horizontal derivative $d_b$. Again by the inner product, we also have
the horizontal gradient $\nabla_b$ with values in $H_g$ and $\Delta_b = \nabla_b^* \nabla_b$.
Note also that dual is taken with respect to integration over $S$, which
means that we also write
$$\int_S |\nabla_b f(\xi)|^2 d\xi = \int_S (\Delta_b f(\xi))\overline{f(\xi)}d\xi,$$
where the norm on the left-hand side is taken in the horizontal
tangent space.

Now the Cayley transform is defined by $(z_0,z) \rightarrow (w_0,w)$, where
$$ w_0 = \frac{z_0 - 1}{ z_0 + 1}, \, w = \frac{2 z }{ z_0 + 1}$$
and we apply this to $z_0 = it + |z|^2$ resp.$z$, where $(z,t)$ is an
element in $H^{2n+1}$ and $|z|^2 = \sum_{j=1}^n z_j \overline{z}_j$. 
This gives the stereographic projection (CR-case)  
of $H^{2n+1} \rightarrow S^{2n+1}$ and also a biholomorphic map from the 
Siegel domain $\{ {\rm Re}\, z_0 > |z|^2\} \rightarrow \{|w_o|^2 + |w|^2 < 1 \}$,
i.e., the complex unit ball. 

On the boundary sphere we get the CR-structure with
contact form $\theta = \frac{i}{2} \sum_{j = 0}^{n+1} w_j d\overline{w}_j$;
the horizontal distribution is given at each tangent space as the maximal
complex subspace, and the metric is that induced from the ambient Euclidian
space; this will be the normalization to be used. The horizontal gradient
is then nothing but the usual gradient of a function, projected orthogonally
onto the horizontal space. On $S = S^{2n+1}$ we again have the CR-Laplacian,
and, adding a suitable constant, this is an intertwining operator between two
principal series representations, namely the CR-Yamabe operator, see
\cite{jl} and \cite{jl2}. Below we shall find the spectrum of this operator.

In the quaternionic case, and also the octonionic case, we have in a similar way
both a noncompact picture on $\overline{N}$ and a compact picture on
$S = K/M$ of a distribution coming from the first summand
in $\mathfrak n = \mathfrak n_{\alpha} \oplus
\mathfrak n_{2\alpha}$. On the sphere this is $K$-invariant and the horizontal
gradient $\nabla_b$ is again obtained via the Euclidean orthogonal projection.
Again, for a suitable constant $\kappa$, $\nabla_b^* \nabla_b + \kappa$ is an
intertwining operator between two principal series representations.                              
          
The first three cases are sometimes called the classical ones;
here we recall the key calculations from \cite{jw} for the eigenvalues of
a $G$-morphism $A_{\nu}$ (intertwining operator) from $Y_{\nu}$ to its dual, normalized
to be $1$ on the constant functions:
\\ \vspace{2 mm}\ \\ 
{\bf The real case.} On spherical harmonics of degree $k$ the eigenvalue is
$$a_k(\nu) = \prod_{j=1}^k \frac{n-1-\nu+j}{\nu+j-1}\,,$$
and for $\nu = (n-2)/2$ this is proportional to $(k+\frac{n-2}{2})(k+\frac{n}{2})$
which ar exactly the expected eigenvalues of the Yamabe operator $\Delta +
\frac{n(n-2)}{4}$; in particular we find the well-known spectrum $k(k+n-1)$ for
the Laplace operator.
\\ \vspace{2 mm}\ \\ 
{\bf The complex case.} The spherical harmonics decompose under $K$ into
${\cal H}_k = \sum_{p+q = k} {\cal H}^{p,q}$ corresponding to holomorphic type
$p$ and anti-holomorphic type $q$. The eigenvalues are
$$ a_{p,q}(\nu) = \prod_{j=1}^p \frac{2n-\nu+2j}{\nu+2j-2} \prod_{l=1}^q \frac{2n-\nu+2l}{\nu+2l-2}\,.$$
We want to consider the case $\nu = n$ in order to find the eigenvalues of the
second-order differential intertwining operator. This gives $(2p + n)(2q + n),$
and subtracting the constant term we obtain the values $4pq + 2(p+q)n = k(k+2n) - j^2$,
where $k = p+q, \, j = p-q$, for the eigenvalues of the CR-Laplacian $\Delta_b$.
Note that this is consistent with standard calculations, see e.g. \cite{f},
where our $\Delta_b = 2 {\rm Re}\Box_b, \, \Box_b = \overline{\partial}_b^* \overline{\partial}_b$            
in terms of the tangential Cauchy--Riemann complex. 
\\ \vspace{2 mm}\ \\ 
{\bf The quaternionic case.} The spherical harmonics decompose in this case under $K$ into
${\cal H}_k = \sum_q V^{k,q}$ corresponding to certain irreducible
representations $V^{p,q}, \, k=p$; the sum is over $p \geq q \geq 0$ and $p-q$ even. We
set $r = (p-q)/2,\, s = (p+q)/2$. Then the eigenvalues are
$$ a_{p,q}(\nu) = \prod_{j=1}^r \frac{4n+2-\nu+2j}{\nu+2j-4} \prod_{l=1}^s \frac{4n+4-\nu+2l}{\nu+2l-2}\,,$$
and we are again interested in a particular $\nu$, namely $\nu=2n+2$ corrresponding
to the second-order differential intertwining operator. This gives
$(2n+2r)(2n+2+2s);$ subtracting the constant part we get for the
eigenvalues of the CR-Laplacian $\Delta_b$ just $k(k+4n+2) - j(j+2)$,
where $2s=k+j,\,2r = k-j, \, k=p, \, p-q = j$.
\\ \vspace{2 mm}\ \\   
{\bf The octonionic case.} Here we use the calculations for this group in \cite{boo}.
We also refer to \cite{j} for the precise relation between spherical harmonics
and the $K$-types occurring in $L^2(S)$, and for more details on the action
of $K$.  
Again we have the eigenvalues of intertwining operators for the spherical
principal series, now found via the method of spectrum generating operators.
(In fact, the same eigenvalues are found in \cite{j} by the same method as in
\cite{jw}.) 
The spectral function (the eigenvalues of the intertwining operator) is in this case
$$Z = a_{k,j}(r) = \frac{\Gamma(j+k+\frac{11}{2}+\frac{r}{2})\Gamma(\frac{11}{2} - \frac{r}{2})
\Gamma(k+\frac{5}{2} +\frac{r}{2})\Gamma(\frac{5}{2} -\frac{r}{2})}{
\Gamma(j+k+\frac{11}{2}-\frac{r}{2})\Gamma(\frac{11}{2} + \frac{r}{2})
\Gamma(k+\frac{5}{2} -\frac{r}{2})\Gamma(\frac{5}{2} +\frac{r}{2})}\,,$$
where $j,k \in \mathbb N$ label the K-types; $\mathfrak k = {\mathfrak so}(9),$
and we label the representations in the usual way via their highest weight
$(\lambda_1, \lambda_2, \lambda_3) = (k + \frac{1}{2}j, \frac{1}{2}j, \frac{1}{2}j, \frac{1}{2}j).$
Note that $K = {\rm Spin}(9)$ and $M = {\rm Spin}(7)$ with a nonstandard
imbedding. The parameter $r$ here is including the $\rho$-shift; we have the
positive root spaces $\mathfrak g_{\alpha}$ and $\mathfrak g_{2\alpha}$ of
dimensions $8$ and $7$ respectively, hence $\rho = 11$ and the above
$\nu = 11 - r$. For the second-order differential intertwining operator we
have $r = 1$, and the relevant eigenvalues are $4(j+k+5)(k+2) - 40$. With
$N = j + 2k$ this is equal to $N(N+14) - j(j+6)$, where $N$ is exactly the
degree of spherical harmonics on $S$ that we decompose under $K$; see \cite{j}
where the $K$-types are labeled $V^{N,j}$ with our notation for the parameters, and 
$N \geq j \geq 0, \, N-j$ even. Summarizing, we have that the eigenvalues of
the CR-Laplacian $\Delta_b$ in this case are $N(N+14) - j(j+6), \, 0 \leq j \leq N$.                 
    
 \section{Logarithmic Sobolev inequalities for rank-one groups}

We can now state our main result in this paper. 
\begin{theorem}
Let $G$ be a split rank-one group and $S = G/P = K/M$ the
corresponding flag manifold; then with the normalized rotation-invariant
measure $d\xi$ on $S,$ we have for any smooth function $f$ on $S$   
(and we may extend naturally by taking suitable limits of functions)
 \begin{equation}
\int_S |f(\xi)|^2 {\rm log}|f(\xi)|d\xi \leq C \int_S |\nabla_b f(\xi)|^2 d\xi + ||f||_2^2 {\rm log}||f||_2 
 \end{equation}
where $\nabla_b$ denotes the boundary CR-gradient, and $||f||_2$ the usual $L^2$-norm.
In the four cases the constant is: 
\begin{itemize}
\item (real case) $C = \frac{1}{n}, \,  G = SO_o(1,n+1), \, S = S^n,$
\item (complex case) $C = \frac{1}{2n}, \, G = SU(1, n+1), \, S = S^{2n+1},$
\item (quaternionic case) $C = \frac{1}{4n}, \, G = Sp(1, n+1), \, S = S^{4n+3},$
\item (octonionic case) $C = \frac{1}{8}, \, G = F_4, \, S = S^{15}.$ 
\end{itemize} 
\end{theorem} 
{\it Proof.} We shall use the inequality found by Beckner
for the sphere $S:$
\begin{equation}
\int_S |F(\xi)|^2 {\rm log}|F(\xi)| \leq \sum_k k \int_S |Y_k(\xi)|^2, 
\end{equation}
for $F = \sum_k Y_k$ decomposed into spherical harmonics; see \cite{be1}
equation (8). Again
we assume $F$ is normalized in $L^2$ i.e., 
$\int_S |F|^2 = 1$. This comes from
the limit $p = 2$ in the HLS inequality. Now we employ the spectrum of the
operator $B = \nabla_b^* \nabla_b = \Delta_b$ (which in the real case
is just $\Delta$) in our four cases:
\begin{itemize}
\item (real case) on ${\cal H}_k$ we have $B = k(k+n-1)$ and the
estimate $$k \leq \frac{k(k+n-1)}{n}\,.$$
\item (complex case) on ${\cal H}_k$ we have $B = k(k+2n) - r^2, \, 
-k \leq r \leq k$ and the
estimate $$k \leq \frac{k(k+2n) - r^2}{2n}\,.$$
\item (quaternionic case) on ${\cal H}_k$ we have $B = k(k+4n+2) -j(j+2), \, 
0 \leq j \leq k$ and the estimate $$k \leq \frac{k(k+4n+2) -j(j+2)}{4n}\,.$$
\item (octonionic case) on ${\cal H}_k$ we have $B = k(k+14) - j(j+6), \,
0 \leq j \leq k$ and the 
estimate $$k \leq \frac{k(k+14) - j(j+6)}{8}\,.$$
\end{itemize}
These tell us that on each degree $k$ of spherical harmonics we have
$k \leq C \Delta_b$ with $C$ as in the theorem, as required. \hfill QED 

There are many important applications of logarithmic Sobolev inequalities,
such as the above; it could be to the Poisson semigroup, to spectral theory, or as in the following
example to smoothing properties of the corresponding heat semigroup, in analogy
with the contraction properties of the Ornstein--Uhlenbeck semigroup.
 
\begin{corollary}
In each of the four cases we have the contraction estimate
for the norm of the semigroup $(t \geq 0)$, ${\rm exp}(-tB): L^q(S) \rightarrow L^p(S)$,
$$||{\rm exp}(-tB)||_{q,p} \leq 1, \, {\rm for}\, {\rm exp}(-t/C) \leq \sqrt{\frac{q-1}{p-1}}$$ 
where $B = \Delta_b$ and $C$ has the value as in the theorem.
\end{corollary}
{\it Proof.} This follows from \cite{gr2} since our $\Delta_b$ is a Sobolev generator. QED
            
 \section {Inequalities in the noncompact picture}
In this section we shall give a new proof of L. Gross' logarithmic
Sobolev inequality on $\mathbb R^n$ with the Gauss measure, using our main theorem.
The idea is to transfer the logarithmic inequality on the sphere to the
flat space (Euclidean space in the real case, and a nilpotent group in
the CR-case). Consider functions only depending on a fixed number of
variables, and let the number of remaining variables tend to infinity. In this
way the Gauss measure turns up in the real case and we obtain the classical
inequality of L. Gross.
  
 \subsection{Stereographic projection in the real case}
 For this section, it is useful as above to realize the representation $\pi_\lambda$ as 
acting on smooth sections of a line
 bundle over $S$, and then to consider the explicit transform to the noncompact picture.
In group-theoretic terms we use the orbit of $\overline{N} = \theta(N)$ in $G/P$ to
provide coordinates; in this way $\pi_{\lambda}$ is realized in functions on $\overline{N}$
and the logarithmic Sobolev inequality becomes translated into an inequality
on $N$ (identified with $\overline{N}$) equipped with a suitable probability
measure and usual CR-structure. Let us first look at the real case (where
the group $N = \mathbb R^n$, and we are dealing with the usual conformal
structure).

The transition from the compact to the noncompact picture is here given
by the stereographic projection
$$x = \frac{\xi}{1+\xi_{n+1}}$$
where $x = (x_1, x_2, \dots , x_n) \in \mathbb R^n$ and
$\xi = (\xi_1, \xi_2, \dots , \xi_n) \in \mathbb R^n, \, \xi_{n+1} \in \mathbb R, \,
(\xi, \xi_{n+1}) \in S^n$, and the inverse is given by
$$\xi = \frac{2x}{1+|x|^2}\,\,, \,\, \xi_{n+1} = \frac{1-|x|^2}{1+|x|^2}\,\,.$$
This is conformal with conformal factor $1+\xi_{n+1} = \frac{2}{1+|x|^2}$.
Hence the Euclidean measures on $S^n$ resp. $\mathbb R^n$
 are related by $d\xi  = \left(\frac{2}{1+|x|^2}\right)^{-n} dx$
and the norm of the gradients will scale in a similar way: Since by definition
$|dF|^2 = |\nabla F|^2$, and since the inner product in the cotangent space
scales with $\lambda^{-2}$ when the inner product in the tangent space scales
with $\lambda^2$, we obtain the following form in $\mathbb R^n$ of the
logarithmic Sobolev inequality on the sphere: 
$$c_n \int_{\mathbb R^n} |f(x)|^2 {\rm log}|f(x)| (1+|x|^2)^{-n} dx \leq
\frac{c_n}{4n} \int_{\mathbb R^n} |\nabla f(x)|^2 (1+|x|^2)^{-n+2} dx$$
for $c_n\int_{\mathbb R^n} |f(x)|^2 (1+|x|^2)^{-n} dx = 1$ and
$c_n\int_{\mathbb R^n} (1+|x|^2)^{-n} dx = 1$.

Now we consider the change of variable to $x/\sqrt{n}$ and functions              
of the form $x \rightarrow f(x/\sqrt{n})$, using $\nabla (f(x/\sqrt{n})) =
\frac{1}{\sqrt{n}}(\nabla f)(x/\sqrt{n})$; this gives
$$c'_n \int_{\mathbb R^n} |f(x)|^2 {\rm log}|f(x)| \left(1+\frac{|x|^2}{n}\right)^{-n} dx $$ $$\leq
\frac{c'_n}{4} \int_{\mathbb R^n} |\nabla f(x)|^2 \left(1+\frac{|x|^2}{n}\right)^{-n+2} dx$$
for $c'_n\int_{\mathbb R^n} |f(x)|^2 \left(1+\frac{|x|^2}{n}\right)^{-n} dx = 1$ and
$c'_n\int_{\mathbb R^n} \left(1+\frac{|x|^2}{n}\right)^{-n} dx = 1$. 

In order to evaluate the normalization constants needed here and later
we record the following. 

\begin{lemma}
For $2N > m$ we have
$$\int_{\mathbb R^m} \left(1 + \frac{|x|^2 + |y|^2}{n} \right)^{-N} dy = 
n^{m/2} \pi^{m/2} \frac{\Gamma(N - \frac{m}{2})}{\Gamma(N)} \left(1+ \frac{|x|^2}{n}\right)^{-N+\frac{m}{2}}$$
for any $|x|^2 \geq 0$.  
\end{lemma}   
As a next step we fix $k$ and let $n = m+k$ be large; assume the function $f(x)$
only depends on the first $k$ variables: $f(x) = f(x_1, x_2, \dots , x_k)$ so that
we can perform the integration in the remaining variables first. With the notation
$x \in \mathbb R^k, \, y \in \mathbb R^m,$ we calculate
\begin{align*}
I &= \int_{\mathbb R^m} \left(1+\frac{|x|^2 + |y|^2}{n}\right)^{-n} dy\\
  &= \left(1+\frac{|x|^2}{n}\right)^{-n} \int_{\mathbb R^m} \left( 1 + \frac{|y|^2/n}
{1 + \frac{|x|^2}{n}}\right)^{-n} dy
\end{align*}
where we change variables to $y/\sqrt{(n(1+\frac{|x|^2}{n})}$ in order to get
$$I = d_{n,k} \left(1+\frac{|x|^2}{n}\right)^{(-n-k)/2}$$        
with the normalizing constant satisfying
$$d'_{n,k} \int_{\mathbb R^k} \left(1 + \frac{|x|^2}{n}\right)^{(-n-k)/2} dx = 1$$
and $d'_{n,k} = c'_n d_{n,k}$.
In fact, from the lemma we find
$$c'_n = n^{-n/2} \pi^{-n/2} \frac{\Gamma(n)}{\Gamma(\frac{n}{2})}$$
$$d_{n,k} = n^{m/2}\pi^{m/2}\frac{\Gamma(\frac{n+k}{2})}{\Gamma(n)}$$
$$ d'_{n,k} = n^{-k/2} \pi^{-k/2}\frac{\Gamma(\frac{n+k}{2})}{\Gamma(\frac{n}{2})}$$
and we shall also use below
$$\tilde{d}_{n,k} = n^{m/2}\pi^{m/2} \frac{\Gamma(\frac{n+k}{2} -2)}{\Gamma(n-2)}$$
as well as $\tilde{d}'_{n,k} = c'_n \tilde{d}_{n,k}$.
      
Integrating the $y$-variable in our inequality we obtain
$$d'_{n,k} \int_{\mathbb R^k} |f(x)|^2 {\rm log}|f(x)| \left(1+\frac{|x|^2}{n}\right)^{(-n-k)/2} dx$$
$$ \leq
\frac{\tilde{d}'_{n,k}}{4} \int_{\mathbb R^k} |\nabla f(x)|^2 \left(1+\frac{|x|^2}{n}\right)^{(-n-k)/2+2} dx$$
for 
$$d'_{n,k}\int_{\mathbb R^k} |f(x)|^2 \left(1+\frac{|x|^2}{n}\right)^{(-n-k)/2} dx = 1$$ 
and (again)  
$$d'_{n,k}\int_{\mathbb R^k} \left(1+\frac{|x|^2}{n}\right)^{(-n-k)/2} dx = 1.$$
Now we take the limit $n \rightarrow \infty$, taking into account the asymptotics
of $d'_{n,k} \sim (2\pi)^{-k/2}$ and 
$\tilde{d}'_{n,k}/4 \sim (2\pi)^{-k/2}$
from $\frac{\Gamma(z+a)}{\Gamma(z)} \sim z^a, \, |z| \rightarrow
\infty$, and also $(1 + \frac{a}{n})^{-n} \rightarrow e^{-a}$, and we finally obtain
$$\int_{\mathbb R^k} |f(x)|^2 {\rm log}|f(x)| d\nu(x) \leq \int_{\mathbb R^k} |\nabla f(x)|^2 d\nu(x)$$
for the Gauss measure $d\nu(x) = (2\pi)^{-k/2} e^{-|x|^2/2} dx$ and $\int_{\mathbb R^k} |f(x)|^2 d\nu(x) = 1$.
Hence we have arrived at the logarithmic Sobolev inequality of L. Gross, as in the appendix.

\subsection{Cayley transform in the usual CR-case} 
Here we employ the Cayley transform between the Heisenberg group and the
CR-sphere; as in the real case we transform the horizontal gradient and the
measure to the Heisenberg group, and get the corresponding form of the
logarithmic Sobolev inequality. Recall the explicit coordinate changes 
$$w_0 = (z_0 - 1)/(z_0 + 1), \, w = 2z/(z_0 + 1)$$ where $z_0 = it + |z|^2$
so that the measure, see e.g., \cite{jl2}, on $S = S^{2n+1}$ becomes (up to a constant)
$((1+|z|^2)^2 + t^2)^{-n-1} dz dt$ with $dz$ and $dt$ denoting Lebesgue measures
in $\mathbb C^n$ resp. $\mathbb R$. 

On the Heisenberg group we have the left-invariant CR-holomorphic vector
fields
$$ Z_j = \frac{\partial}{\partial z_j} + i \overline{z_j} \frac{\partial}{\partial t}$$
corresponding to the distribution; the real and imaginary parts form a basis
of the distribution and define the CR-gradient. 
Explicitly we have the real left-invariant vector fields
$$ X_j = \frac{\partial}{\partial x_j} + 2 y_j \frac{\partial}{\partial t}, \,
   Y_j = \frac{\partial}{\partial y_j} - 2 x_j \frac{\partial}{\partial t}$$
for $j = 1, 2, \dots , n$. Then the CR-Laplacian is also
$\Delta_b = \sum_1^n (X_j^2 + Y_j^2)$.
   
We can make the change of
variables similar to the real situation as $f(\frac{z}{\sqrt{2n}}, \frac{t}{2n})$; now
$$\nabla_b f\left(\frac{z}{\sqrt{2n}}, \frac{t}{2n}\right) = \frac{1}{\sqrt{2n}}(\nabla_b f)
\left(\frac{z}{\sqrt{2n}}, \frac{t}{2n}\right)$$ 
and furthermore  again
we have to take into account how the CR-gradient changes by the CR-conformal
factor, see \cite{jl2}. In this way we can 
write the logarithmic Sobolev inequality in the CR-case on the Heisenberg group
as
$$c'_n\int_{\mathbb C^n \times \mathbb R} |f(z,t)|^2 {\rm log}|f(z,t)| d\mu_n(z,t) \leq
\frac{c'_n}{4}\int_{\mathbb C^n \times \mathbb R} |\nabla_bf(z,t)|^2 d\mu_{n-1}(z,t)$$
for $c'_n\int_{\mathbb C^n \times \mathbb R} |f(z,t)|^2 d\mu_n(z,t) = 1$. The
measure is here the (up to the constant $c'_n$ probability) measure
$$d\mu_n(z,t) = \left(\left(1+\frac{|z|^2}{2n}\right)^2 + \frac{t^2}{4n^2}\right)^{-n-1}dz dt$$ 
with $dz$ and $dt$ Lebesgue measures as before; here
$$c'_n \int_{\mathbb C^n \times \mathbb R} d\mu_n(z,t) = 1$$
$$c'_n = (2n)^{-n-1} \pi^{-n-\frac{1}{2}}\frac{\Gamma(2n+1)}{\Gamma(n + \frac{1}{2})}$$
by standard tables, e.g., \cite{gr} p. 343 formula 2.    
This provides a new inequality on the Heisenberg group, and it might be possible
to obtain some analogue of the Gaussian logarithmic Sobolev inequality on $\mathbb R^k$ as a
consequence. We shall refrain from completing this idea, but just limit ourselves
to giving a few explicit inequalities that one may immediately deduce.          
 
Now in order to see what happens, we try the trick that
we used in the real case, namely that of letting the function
only depend on the first $k$ variables. Thus we will write $z = (u,v) \in \mathbb C^k \times
\mathbb C^m, \, n = m+k$ with $k$ fixed and $n$ large, and consider the integral
$$I = \int_{\mathbb C^m} \left( \left(1 + \frac{|u|^2+ |v|^2}{2n}\right)^2 + \frac{t^2}{4n^2}\right)^{-N} dv$$
which we evaluate using \cite{gr} p. 345, formula 10.    
The result is   
\begin{align*}
I &= (2n)^m (1+\frac{|u|^2}{2n})^{-2N+m} \frac{2\pi^m}{\Gamma(m)}\\
&\times\int_0^{\infty}
\left((1+r^2)^2+\frac{t^2}{4n^2}\left(1+\frac{|u|^2}{2n}\right)^{-2}\right)^{-N} r^{2m-1} dr
\end{align*} 
where the last integral for $t=0$ equals $\frac{1}{2} B(m,2N-m) = \frac{\Gamma(m) \Gamma(2N-m)}
{2\Gamma(2N)}$ in terms of the
usual beta function, and in general can be further rewritten as
$$\frac{1}{2}(1+D)^{-N + \frac{m}{2}}\int_0^{\infty} (x^2+2\beta x + 1)^{-N}x^{m-1} dx$$   
where $D = C/A^2, \, C = t^2/(2n)^2,\, A = 1 + \frac{|u|^2}{2n},\, \beta = (1+D)^{-1/2}.$   

Now we integrate with respect to the $v$-variable in the inequality on the Heisenberg group
and find the asymptotics for large $n = m+k$ with $k$ fixed. Summarizing, we obtain
$$\frac{1}{\sqrt{n}}\int_{\mathbb C^k \times \mathbb R} |f|^2 {\rm log}|f| d\nu_n \leq
\frac{1}{\sqrt{n}}\int_{\mathbb C^k \times \mathbb R} |\nabla_bf|^2 d\rho_n
+ 8\sqrt{n} \int_{\mathbb C^k \times \mathbb R} |\frac{\partial f}{\partial t}|^2 \widetilde{d\rho_n}$$
for $\frac{1}{\sqrt{n}}\int_{\mathbb C^k \times \mathbb R} |f|^2 d\nu_n = 1$. Here
$\frac{1}{\sqrt{n}}\int_{\mathbb C^k \times \mathbb R} d\nu_n = 1$ and we have the
asymptotic relations
$$ d\rho_n(u,t) \sim \widetilde{d\rho_n}(u,t) \sim 
   d\nu_n(u,t) \sim (2\pi)^{-k} e^{-|u|^2/2} \frac{du\, dt}{\sqrt{2\pi}}$$
for $n \rightarrow \infty$.
     
 \section*{Appendix : Hardy--Littlewood--Sobolev inequalities}

These are inequalities of the following classical type \cite{hlp}: Suppose we have
two finite sequences of nonnegative real numbers $a_i, b_i, i = 1, \dots ,n,$
and we consider the sum
$$Q = \sum_{i = 1}^n a_i b_i.$$
Then with the same sequences rearranged in decreasing order,
$a_1^* \geq a_2^* \dots \geq a_n^*,$ resp. $b_1^* \geq b_2^* \dots \geq b_n^*,$ 
we have that $Q^* \geq Q$ where now
$$Q^* = \sum_{i = 1}^n a_i^* b_i^*.$$
As we see in \cite{hlp} the same principle of rearrangement may be
extended to functions and many other types of expressions; similar
ideas are found in other forms of symmetrization, such as e.g., 
Steiner symmetrization. We shall be interested in quantities of the form
$$ Q = \int \int f(x)g(y)h(x-y) dx dy$$
with each integration being over $\mathbb R^n$ and $dx, dy$ denote
Lebesgue measure.

Then for nonnegative measurable functions $f, g, h$ we consider their
equimeasurable symmetric non-increasing rearrangements $f^*, g^*, h^*$
(as in \cite{hlp}) and have $Q^* \geq Q$ where now
$$ Q^* = \int \int f^*(x)g^*(y)h^*(x-y) dx dy$$
\cite{hlp}. This forms the basis of E. Lieb's \cite{l} deep analysis
where he establishes the following sharp Hardy--Littlewood--Sobolev 
(HLS) inequality (see also \cite{be1}):  

\begin{proposition}
Let $S = S^n$ be the $n$-dimensional sphere with the normalized
usual rotation-invariant measure, $0 < \lambda < n, \, 0 < p = \frac{2n}{2n - \lambda} < 2;$
then we have the estimate for the $L^p$-norm $||F||_p = (\int_S |F|^p)^{1/p}$ 
\medskip
$$\sum_{k = 0}^{\infty} \gamma_k \int_S |Y_k|^2 \leq ||F||_p^2,$$
where $F = \sum_{k = 0}^{\infty} Y_k$ is the decomposition of the
measurable function $F$ into spherical harmonics $Y_k \in {\cal H}_k$ of degree $k$,
$$\gamma_ k = \frac{\Gamma(\frac{n}{p}) \Gamma(\frac{n}{p'} + k)}
{\Gamma(\frac{n}{p'}) \Gamma(\frac{n}{p} + k)}$$
and $p, p'$ are dual exponents: $\frac{1}{p} + \frac{1}{p'} = 1$.
 
This is equivalent to giving the best constant 
$$K_p = \pi^{n/p'} \frac{\Gamma(\frac{n}{p} - \frac{n}{2})}{\Gamma(\frac{n}{p})} 
\left(\frac{\Gamma(\frac{n}{2})}{\Gamma(n)}\right)^{(p-2)/p}$$ 
in the estimate
$$\int \int f(x)|x - y|^{-\lambda} g(y) dx dy \leq K_p ||f||_p ||g||_p$$
for nonnegative functions and their $L^p$ norms,
where the measures $dx, dy$ are Lebesgue measure and the integrals over
$\mathbb R^n$.
\end{proposition}
Note that the integral $If(y) = \int f(x) |x - y|^{-\lambda} dx$ defines an
intertwining operator between two principal series representations
of $G = SO(1,n+1)$,
namely from a representation to its natural dual.
It is a very important part of the theory that one can find the eigenvalues
on the sphere, that is, in the compact picture of
the principal series representations of this family of intertwining operators. 
Note that when we make a change of variables and transform $I$ to the sphere and normalize it
so that $I1 = 1$, the sharp HLS inequality states that $||If||_{p'}
\leq ||f||_p$. It is an appealing conjecture, that such a contraction
property holds more generally, i.e., for other groups and their
continuations of principal series to suitable real parameters. 
The representations here belong to the complementary series and
they are unitary through the invariant Hermitian form coming from
the intertwining operator. On the other hand, the $L^p(S)$ Banach
norm is invariant in the original space, and the $L^{p'}(S)$
Banach norm is invariant in the target space. So another way to think
of the sharp HLS inequality is the following for the invariant
unitary norm $||f||$ (say for real functions):
$$||f||^2 = <If,f> \leq ||If||_{p'} ||f||_p \leq ||f||_p^2.$$ 
The conjecture would be, that this remains true in the CR-case, where
the complex case is already interesting.          
     
Now as demonstrated in \cite{be1} it is very interesting to study  
the parameter endpoints $p = 1, 2$ in HLS, where one may consider the
derivatives in the parameter. One result that follows at $p=2$ is the celebrated
logarithmic Sobolev inequality below \cite{gr1} for the Gauss measure
$d\nu = (2\pi)^{-n/2} e^{-|x|^2/2} dx$ on $\mathbb R^n.$ There are
several different proofs of this result, and in this paper we have given
a new way of deriving it from the real case in our main Theorem.      
At $p=1$ one obtains \cite{be1} an exponential-class inequality of
Moser-Trudinger type. It is a highly interesting problem to find the right
analogues of exponential-class inequalities in the 
framework of the CR-geometries considered in this paper.    
\begin{proposition} For a smooth function $f$ on $\mathbb R^n$  
(or suitable limit functions)     we have the estimate 
$$\int |f|^2 {\rm log}|f| d\nu \leq \int |\nabla f|^2 d\nu $$
for $\int |f|^2 d\nu = 1$.
\end{proposition}  
More generally, if we have a probability space $(\Omega, \mu)$
and a self-adjoint linear operator $B$ on $L^2(\mu)$ with $B \geq 0$
satisfying
$$\int_{\Omega} |f|^2 {\rm log}|f| d\mu \leq (Bf,f)$$
for all $f$ in the domain of $B$ with $||f||_2 = 1$, then we call this a
logarithmic Sobolev inequality with {\it Sobolev generator} $B$.

 \bigskip
 \footnotesize{ \noindent Address: 
 Matematisk Institut, Byg.\,430, Ny Munkegade, 8000 Aarhus C,
 Denmark.
 \medskip

 \noindent \texttt{
  orsted@imf.au.dk
 }}


\begin{thebibliography}{99}\itemsep=-.2pc

\bibitem{be1}
Beckner W., Sharp Sobolev inequalities on the sphere and the Moser-trudinger
inequality, {\it Ann. of Math. (2)}, {\bf 138} (1993), no. 1, 213--242. 
 
\bibitem{be2} Beckner W., Geometric inequalities in Fourier analysis,
{\it Essays on Fourier analysis in honor of Elias M. Stein} (Princeton, NJ, 1991),
36--68, Princeton Math. Ser., 42, Princeton University Press, Princeton, NJ, 1995. 

\bibitem{boo} Branson T., \'O{}lafsson G., and \O{}rsted B., Spectrum generating
operators and intertwining operators for representations induced from a
maximal parabolic subgroup, {\it Journ. Func. Anal.} {\bf 135} (1996), no. 1,
163--205. 

\bibitem{f} Folland G. B., The tangential Cauchy-Riemann complex on spheres,
{\it Trans. AMS} {\bf 171} (1972), 83--133. 

\bibitem{gr} Gradsheteyn I. S. and Ryzhik I. M., {\it Table of Integrals, Series, and
Products}, fifth edition, Academic Press, 1994. 
                
\bibitem{gr1} Gross L., Logartihmic Sobolev inequalities, {\it Amer. J. Math.} {\bf 97}
(1975),
 1061--1083. 

\bibitem{gr2} Gross L., Logarithmic Sobolev inequalities and contractivity
properties of semigroups, C.I.M.E. 1992, L.N.M. {\bf 1563}, Springer, 1992.   

\bibitem{hlp} Hardy G. H., Littlewood J. E. and P\'o{}lya G, {\it Inequalties},
Cambridge University Press, 1952.  
 
 \bibitem{jl}
Jerison D. and Lee J. M., Extremals for the Sobolev inequality on the
Heisenberg group and the CR Yamabe problem, {\it J. Amer. Math. Soc.} {\bf 1}
(1988), no. 1, 1--13. 

\bibitem{jl2} Jerison D. and Lee J. M., The Yamabe problem on CR manifolds,
{\it J. Differential Geom.} {\bf 25} (1987), no. 2, 167--197. 
 
\bibitem{j} Johnson K. D., Composition series and intertwining operators
for the spherical principal series, II, {\it Trans. AMS} {\bf 215} (1976),
269--283. 
  
\bibitem{jw} Johnson K. D. and Wallach N., Composition series and intertwining
operators for the spherical principal series, I, {\it Trans. AMS} {\bf 229}
(1977), 137--173. 
 
\bibitem{l} Lieb E. H., Sharp constants in the Hardy-Littlewood-Sobolev and
related inequalities, {\it Ann. of Math.} {\bf 118} (1983), 349--374. 
  
\bibitem{wal} Wallach N., {\it Harmonic Analysis on Homogeneous Spaces}, Marcel Dekker, 1972.

\bibitem{w} Warner G., {\it Harmonic Analysis on Semi-simple Lie Groups I}, Springer Verlag, 1972.

 \end{thebibliography}
 \end{document}